\documentclass[a4paper,10pt]{article}
\usepackage[english]{babel}
\usepackage[latin1]{inputenc}
\usepackage[T1]{fontenc}
\usepackage{amsmath}
\usepackage{amsfonts}
\usepackage{dsfont}
\usepackage{amsthm}
\usepackage{amssymb}
\title{Exponential decay of correlations for generic birational maps of $\P^k$}
\usepackage{dsfont}
\author{Gabriel Vigny}
\begin{document}
\newtheorem{Theorem}{Theorem}
\newtheorem*{Theorem*}{Theorem}
\newtheorem{theorem}{Theorem}[section]
\newtheorem{proposition}[theorem]{Proposition}
\newtheorem{defi}[theorem]{Definition}
\newtheorem{corollaire}[theorem]{Corollary}
\newtheorem{Hypothesis}[theorem]{Hypothesis}
\newtheorem{lemme}[theorem]{Lemma}
\newtheorem{Remark}[theorem]{Remark}
\newcommand{\U}{\mathcal{U}}
\newcommand{\C}{\mathcal{C}}
\renewcommand{\P}{\mathbb{P}}
\renewcommand{\phi}{\varphi}
\newcommand{\Cc}{\mathbb{C}}
\newcommand{\Nn}{\mathbb{N}}
\newcommand{\Rr}{\mathbb{R}}
\newcommand{\Qq}{\mathbb{Q}}
\newcommand{\Zz}{\mathbb{Z}}
\newcommand{\Acal}{\mathcal{A}}
\newcommand{\Bcal}{\mathcal{B}}
\newcommand{\Dcal}{\mathcal{D}}
\newcommand{\Ecal}{\mathcal{E}}
\newcommand{\Gcal}{\mathcal{G}}
\newcommand{\Hcal}{\mathcal{H}}
\newcommand{\Lcal}{\mathcal{L}}
\newcommand{\Ical}{\mathcal{I}}
\newcommand{\Pcal}{\mathcal{P}}
\newcommand{\Qcal}{\mathcal{Q}}
\newcommand{\Scal}{\mathcal{S}}
\newcommand{\Zcal}{\mathcal{Z}}
\newcommand{\supp}{\mathrm{Supp}}
\hyphenation{plu-ri-sub-har-mo-nic}
\date{}

\maketitle
\begin{abstract} We prove the exponential decay of correlations for $C^\alpha$-observables ($0<\alpha \leq 2$) for generic birational maps of $\P^k$ \emph{\`a la} Bedford-Diller. In the particular case of regular birational maps, we give a better estimate of the speed of the decay, getting results as sharp as Dinh's results for H\'enon maps. 
\end{abstract}


\section{Introduction}
Let $f$ be a (dominant) generic birational map of $\P^k(\Cc)$ \emph{\`a la} Bedford-Diller of algebraic degree $d>1$. In dimension $2$, those maps were introduced by Bedford and Diller in \cite{BD1}. A map $f$ is in that set of maps if it satisfies the condition:
\begin{align*}
&\sum^{\infty}_{n=0}      \frac{1}{d^n}\log \text{dist}(I(f),f^n(I(f^{-1}))) >-\infty \\
& \qquad \qquad \qquad \mbox{and} \\
&\sum^{\infty}_{n=0}  \frac{1}{d^n}\log \text{dist}(I(f^{-1}),f^{-n}(I(f))) >-\infty,
\end{align*}
where $I(f^{\pm 1})$ denotes the indeterminacy set of $f^{\pm 1}$. For such a map, the authors constructed an invariant probability measure $\mu$ as the intersection of the Green currents $T^+_1$ and $T^-_1$. The currents $T^\pm_1$ are defined as $\lim_{n\to \infty} d^{-n} (f^{\pm n})^*(\omega)$ where $\omega$ is the Fubini-Study form on $\P^2$. The measure $\mu$ is mixing \cite{BD1}, hyperbolic and of maximal entropy \cite{Duj}. The above condition is generic in the sense that starting from any birational map $f$ of $\P^2$ of algebraic degree $d \geq 2$, then outside a pluripolar set of $A\in \mathrm{Aut}(\P^2)$, the map $f\circ A$ satisfies the condition. 

De Th\'elin and the author generalized in \cite{DV1} that family to birational maps of $\P^k$ for which $\text{dim}(I(f))=k-s-1$ and $\text{dim}(I(f^{-1}))=s-1$ for some $1\leq s \leq k-1$ (as indeterminacy sets have at least codimension 2, in dimension 2 the condition is always satisfied with $s=1$). In that case, $f^{-1}$ is of algebraic degree $\delta$ with $d^s=\delta^{k-s}$. We gave a condition similar to the above one for which we constructed an invariant probability measure $\mu$. The Green current $T^+_s$ (resp. $T^-_{k-s}$) of bidegree $(s,s)$ (resp. $(k-s,k-s)$) is defined by $\lim_{n\to \infty} d^{-sn} (f^{n})^*(\omega^s)$ (resp. $\lim_{n\to \infty} d^{-sn} (f^{- n})^*(\omega^{k-s})$) and the measure $\mu$ is $T^+_s \wedge T^-_{k-s}$.
The measure $\mu$ is mixing, hyperbolic and of maximal entropy. The condition we gave is generic in the same sense than in dimension 2 (see also \cite{moi4} for a generalization of that result to arbitrary rational maps). We call generic birational map of $\P^k$ such a map.
An interesting particular case of generic birational map is given by the so-called \emph{regular birational maps} of $\P^k$ introduced by Dinh and Sibony in \cite{DS10}. That family is open (stable under small perturbations by automorphisms) and contains the regular automorphisms of $\Cc^k$ (\cite{Sib}).  \\

In all the above articles, mixing is achieved using the fact that the Green currents are extremal. Although beautiful, the proof is not direct and uses  extraction of subsequences. In particular, no speed of mixing is obtained. Our main result here is the following theorem of decay of correlations for $C^\alpha$-observables ($0<\alpha\leq 2$). Not only does it implies mixing, but it gives precise estimates which is what is important in practice (for example, for meromorphic maps of large topological degree, exponential decay implies the central limit theorem for bounded q.p.s.h-observables  \cite{DS4}). Finally, one of the main topic in dynamics is to classify the different dynamics, namely we want to know if two different dynamics are conjugated and what degree of regularity can be achieved for the conjugacy. Decay of correlations for $C^\alpha$-observables can give evidences that two dynamics are not conjugated in $C^\alpha$, since if one dynamics has exponential decay for such observables and not the other (or with a different exponent), then they cannot be conjugated (mixing only gives results for mesurable conjugacy). We now state our theorem: 
\begin{theorem}\label{decay}
 Let $\varphi$ and $\psi$ be two functions in $C^\alpha$, $0< \alpha \leq 2$. Then there exists a constant $C_\alpha$ independent of $\varphi$, $\psi$ such that for all $N\in \Nn$:
 \begin{equation}\label{eqdecay}
 \Big| \mu(\varphi\circ f^N . \psi)- \mu(\varphi)\mu(\psi)\Big|\leq C_\alpha \| \varphi \|_{C^\alpha} \|\psi\|_{C^\alpha} d^{-\frac{\alpha sN}{4k}} 
 \end{equation}
 Assume furthermore that $f$ is a regular birational map, then we have:
 \begin{equation}\label{eqdecay2}
 \Big| \mu(\varphi\circ f^N . \psi)- \mu(\varphi)\mu(\psi)\Big|\leq C_\alpha \| \varphi \|_{C^\alpha} \|\psi\|_{C^\alpha} d^{-\frac{\alpha sN}{2k}} 
 \end{equation}
 \end{theorem} 
Exponential decay of correlations is a classical question in complex dynamics and some results have been achieved. This paper is particularly interested in Dinh's results (\cite{Dinh}) where he proves exponential decay of correlations for regular automorphisms of $\P^k$. He made the technical hypothesis that $2s=k$ (in particular, $k$ is even, that is the case for H\'enon maps of $\P^2$ where $s=1$). Nevertheless, his proof works for any regular automorphism with little modifications. We now benchmark our results with his. Consider the above theorem:
\begin{itemize}
	\item In the case where $k=2$, then we obtain a decay in $d^{-\frac{N}{4}}$ in the generic case.
	\item In the regular case, assume as in \cite{Dinh}[Theorem 1.1] that $2s =k$, then we get a decay in $d^{-\frac{N}{2}}$ and our results extend the ones of Dinh (and are thus as sharp as them).  
\end{itemize}
Now, it is not clear whereas the decay we have in the generic case is the sharpest one can get. A better statement would be to prove the following:
\begin{equation}\label{eqdecay3}
\Big| \mu(\varphi\circ f^N . \psi)- \mu(\varphi)\mu(\psi)\Big|\leq C_\alpha \| \varphi \|_{C^\alpha} \|\psi\|_{C^\alpha} d^{-\frac{\alpha sN}{2k}} 
\end{equation} 
 for all $0<\alpha \leq 1$ ($\alpha \leq 2$ seems too optimistic). \\
 
Exponential decay has also been proved in the following birational settings: for  automorphisms of compact K\"ahler manifolds (\cite{DS12}) and for invertible horizontal-like maps (\cite{DinhNguyenSibony1}[Theorem 5.1]). Note also the case of meromorphic maps of compact K\"ahler manifolds of large dynamical degree (\cite{DS4}) which has inspired this article. \\
  
 Let us sketch (very roughly) the proof of the decay of correlations in \cite{Dinh}. The author first proves exponential rates of convergence to the Green currents and to the Green measure for the limit of normalized pull-back   and push-forward of suitable currents and their wedge product (using a $dd^c$-argument). He then considers the map $F=(f,f^{-1 })$ which acts a priori on $(\P^k)^2$ but he shows that it can in fact be seen as acting on $\P^{2k}$ as a regular automorphism of $\Cc^{2k}$ (this is where the hypothesis $2s=k$ is used, one has to work in $(\P^k)^2$ otherwise). Write $\Cc^{2k}=\Cc^k \times \Cc^k$ and let $x$ and $y$ be the coordinates on each term of the product. Then he observes:
 \begin{align}\label{letrucdeladiagonale}
\mu(\varphi \circ f^{2n} \psi ) = \mu(\varphi \circ f^n\psi \circ f^{-n})= \int_{\P^{2k}} (\varphi(x) \psi(y)) \circ F^n(x,y) T^+_s(x)\wedge T^-_{k-s}(y) \wedge [\Delta]
\end{align}
where $[\Delta]$ is the diagonal on $\Cc^{2k}$. Using the invariance of the Green currents, he writes it down as $\langle \varphi(x) \psi(y) T^+_s(x) \wedge T^-_{k-s}(y) , (d^{-2sn}F^n)_*[\Delta] \rangle$ and he then uses the fact that $[\Delta]$ is a suitable current (that is the estimates he had for the convergence of the measure apply here).  The strategy is the same in \cite{DS12, DinhNguyenSibony1}. \\

In this article, we cannot use that approach as the map $F$ would be a birational map of $(\P^k)^2$ that does not deal nicely with the diagonal. Instead, we do write a formula similar to (\ref{letrucdeladiagonale}) but we decompose the diagonal into two terms: its cohomological part plus a term in $dd^cV$ where  $V$ is a quasi-potential of the diagonal. The first term is controlled by the results of Section \ref{Construction of the measure} where we give a new proof of the existence of the measure $\mu$ with an explicit rate of convergence. More precisely, we show that for a test function in $C^1$ the sequence
$\langle \varphi \circ f^n T^+_s, \omega^{k-s}\rangle$ converges to $\langle \mu, \varphi \rangle$ exponentially fast (see Proposition \ref{rate_measure}), we give a better estimate in the regular case for $C^2$-observables. The ideas are inspired by the $d$-method of \cite{DS4} but, here, we work on the current $T^+_s$ instead than in the ambient $\P^k$. For the second term, we use suitable Stokes formulae and Cauchy-Schwarz inequalities in Section \ref{Decayofcorr} in order to bound the term by quantity of the type:
$$\int_{(\P^k)^2}  (f^n)^*(\omega(x) \wedge  T^+_s(x)) \wedge T^-_{k-s}(y) \wedge V $$
which can be interpreted as $\langle U_{ T^-_{k-s}}, (f^n)^*(\omega(x) \wedge  T^+_s(x)) \rangle$ where $U_{ T^-_{k-s}}$ is the Green quasi-potential of  $T^-_{k-s}$. \\

The main difficulties in this article lie in the fact that we are dealing with currents of higher bidegree (and not just $(1,1)$ currents) and the indeterminacy sets are not trapped in some particular region of the space (it is the case for regular maps though) and can very well be on the support of the measure. In particular, most of the manipulations (wedge product, Stokes formula \dots) we want to do are not clearly defined. In order to overcome that difficulty, we regularize almost everything and use special dynamical cut-off functions. We then use pluri-potential and super-potential (\cite{DS6}) theories to pass to the limit. In  Section \ref{setting}, we start by recalling the facts we need on super-potential theory and on generic and regular birational maps of $\P^k$. We also give the construction of the cut-off functions we use. The reader unfamiliar with super-potential theory can stick to the case $k=2$ where all the currents are of bidegree $(1,1)$.

\section{Settings}\label{setting}
We recall basics on super-potentials (\cite{DS6} and the appendix of \cite{DV1}). Let $\C_q$ denote the set of positive closed currents of mass $1$ and bidegree $(q,q)$ in $\P^k$ for $0\leq q \leq k$. For $T \in \C_q$, we consider some quasi-potential $U_T$ of $T$ (that is $T=\omega^q+dd^cU_T$, see \cite{DS6}[Theorem 2.3.1]). We define the super-potential $\U_T$ associated to the quasi-potential $U_T$  as the function on $\C_{k-q+1}$ defined for $S$ smooth by: 
$$\U_T(S):=\langle U_T, S\rangle.$$
This definition can be extended by sub-harmonicity along the structural varieties to any $S \in \C_{k-q+1}$ allowing the value $-\infty$. The mean of a super-potential is then $m_T:=\U_T(\omega^{k-q+1})$ and the super-potential of a current is uniquely defined by its mean. This is one of the strengths of super-potentials: quasi-potentials differ by a $dd^c$-closed form which can be a highly irregular object whereas super-potentials are defined up to a constant. We say that a sequence $T_n \in \C_q$ \emph{converges to $T$ in the Hartogs' sense}  (or H-converges) if $T_n \to T$ in the sense of currents and if we can choose super-potentials $(\U_{T_n})$ and $\U_{T}$ such that $m_{T_n} \to m_T$ and $\U_{T_n} \geq \U_T$ for all $n$. In that case, if $S_n \to S$ in the Hartogs' sense in $\C_{k-q+1}$, then $\U_{T_n}(S_n) \to  \U_{T}(S)$. 
 A current is said to be more \emph{ H-regular} than another one if for the right means, its super-potential is greater than the other. All the classical tools (intersection, pull-back), well defined for smooth forms can be extended as continuous objects for the Hartogs' convergence and any such operation well defined for a given current $T$ is also well defined for a more H-regular current:
\begin{itemize} 
\item We say that $T\in \C_q$ and $S\in \C_{k-q-r}$ are wedgeable if $\U_{T}(S\wedge \Omega) >-\infty$ for some smooth $\Omega \in \C_{r+1}$. That condition is symmetric in $T$ and $S$ and if it is satisfied we can define the wedge product $T \wedge S$ with the above continuity property.
If $T'$ and $S'$ are more H-regular than $T$ and $S$ then $T'$ and $S'$ are wedgeable and $T' \wedge S'$ is more H-regular than $T \wedge S$. 
If $R\in \C_r$, $S\in \C_s$ and $T\in \C_t$ ($r+s+t \leq k$) are such that $R$ and $S$ are wedgeable and $R\wedge S$ and $T$ are wedgeable then the wedge product $R \wedge S \wedge T$ is well defined and that property is symmetric in $R$, $S$ and $T$. 
\item Similarly, for $T \in \C_q$, we say that $T$ is $f^*$-admissible if its super-potential is finite at some current of the form $d_{q+1}^{-1}f_*(S)$ for $S\in \C_{k-q+1}$ smooth near $I(f)$ ($d_{q+1}$ is the normalization so that  $d_{q+1}^{-1}f_*(S)$ is of mass $1$). For such current, we can define its pull-back with the above continuity property. If $T'$ is more H-regular than $T$ then $T'$ is also  $f^*$-admissible and  $d_q^{-1}f^*(T')$ is more H-regular than $d_q^{-1}f^*(T)$.
\end{itemize}

In \cite{DV1}, intersections of currents are defined in the sense of super-potentials. It turns out that every object can be written in terms of intersection of $(1,1)$ currents and quasi-plurisubharmonic functions (qpsh for short). 
In that setting, if $T \in \C_q$ and $S=\omega+dd^c u \in \C_1$, then the wedge product $S\wedge T$  is defined by $\omega \wedge T +dd^c(uT)$ provided that $u \in L^1(T \wedge \omega^{k-s})$. In particular, if $S_i=\omega+dd^c u_i \in \C_1$ for $i\leq j$, the wedge product $S_1\wedge S_2 \wedge \dots \wedge S_j$ is defined inductively provided that $u_i\in L^1(S_{i+1} \wedge \dots \wedge S_j \omega^{k+1-j-i})$. It turns out that the wedge product is symmetric with respect to the $S_i$ and $u_{i} \in L^1(S_{i_1}\wedge\dots \wedge  S_{i_m} \wedge \omega^{k-m})$ for all choices of $i_l$. This definition coincides with the definition by the super-potentials so this wedge product depends continuously on the $S_i$ for the Hartogs' convergence(\cite{DS6}[Section 4.4]). Note also that Hartogs' convergence of $S_{i,m}$ to $S_i$ in that setting is implied by the fact  we can choose the quasi-potential $u_{i,m}$ of $S_{i,m}$ such that it decreases to the quasi-potential of $S_i$. Finally, we have that $S'$ is more $H$-regular than $S$ if and only if we can choose quasi-potentials $u'$ of $S'$ and $u$ of $S$ such that $u'\geq u$. \\

We now give the properties of the generic birational map $f$ we will need (see \cite{DV1}[Corollary 3.4.11, Theorem 3.4.1, Theorem 3.4.13]). 
As we mentioned in the introduction, we have that $f^{-1}$ is of algebraic degree $\delta$ with $d^s=\delta^{k-s}$. For $0\leq q \leq k$, there exists $d_q$ such that for all $S \in \C_q$ $f^*$-admissible then $f^*(S)$ is of mass $d_q$.  
For all $q\leq s$ and all $S \in \C_q$ $f^*$-admissible then $f^*(S)$ is of mass $d^q$, similarly for all $q \leq k-s$ and  all $S \in \C_q$ $f_*$-admissible (that is to say $(f^{-1})^*$-admissible) then $f_*(S)$ is of mass $\delta^q$. In particular, $d_q =d^q$ for $q\leq s$ and $d_q= \delta^{k-q}$ for $q \geq s$.
The Green current $T^+_1$ is a positive closed current of order 1. It is the limit in the Hartogs' sense of $1/d^n (f^n)^*(\omega)$ and it is invariant ($f^*T^+_1=dT^+_1$) (see \cite{Sib}). Similarly, $T^-_1=\lim_{n\to \infty} \delta^{-n} (f^n)_*(\omega) \in \C_1$ satisfies $f_*T^-_1=\delta T^-_1$. Then, $T^+_s\in \C_s$ the Green current of order $s$ is equal to $(T^+_1)^s$ and is invariant $f^*T^+_s=d^sT^+_s$ (similarly $T^-_{k-s}=(T^-_1)^{k-s} \in \C_{k-s}$ satisfies $f_*T^-_{k-s}=\delta^{k-s}T^-_{k-s}$). Finally, $T^+_s$ and $T^-_{k-s}$ are wedgeable and the intersection $\mu:=T_s^+ \wedge T^-_{k-s}$ is an invariant probability measure. We have that $\mu $ integrates the quasi-potentials of $T^+_1$ and $T^-_1$ (and thus $\log \mathrm{dist}(x, I^+)$ and $\log \mathrm{dist}(x, I^-)$). In the formalism of super-potentials, this means that $\U_{T^+_1}(T^+_s \wedge T^-_{k-s})>-\infty$.  Standard arguments of Hartogs' regularity imply that for any $S\in \C_1$ more H-regular than $T^+_1$ (in particular for $T^+_s$, $\omega$ and the current $T_n$ defined below) then $S$ and $T^+_s$ are wedgeable and $S \wedge T^+_s$ is $(f^n)^*$-admissible and for any $Q$ more H-regular than $T^-_s$ then $\U_S( T^+_s\wedge Q)$ is finite. 
The measure $\mu$ gives no mass to pluripolar sets (hence to proper analytic sets). \\

If, in addition, $f$ is a regular birational map then we also have that the super-potential $\U_{T^+_s}$ (resp. $\U_{T^-_{k-s}}$) is uniformly bounded at any current with support in $\mathrm{supp}(T^-_{k-s})$ (resp. $\mathrm{supp}(T^+_{s})$).  This is in fact a weaker property than the one satisfied by regular maps for which the Green currents $T^+_j$, $j\leq s$, are PC in a neighborhood of $\mathrm{supp}(T^-_{k-s})$ (see \cite{DS10}[Theorem 1.2], the same property holds for the $T^-_j$, $j\leq k-s$). The sharp estimates we have in the regular case extend to such a map (but we do not have any meaningful examples outside regular maps).\\

We will need the following cut-off function $\chi_A$, its construction follows ideas of \cite{Sibony2}. Let $h\in C^\infty(\Rr)$ be such that $0\leq h \leq 1$, $h= 0$ for $x\leq -2$ and $h = 1$ for $x\geq -1$. Let $v_n$ be a quasi-potential of $(1/d)^{n} (f^*)^{n}(\omega)$. Then $v_n$ is smooth outside $I(f^{n})$ and has singularities in $O(\log \mathrm{dist}(x, I(f^{n})))$ (see \cite{DV1}[Lemma 3.2.4]). Substracting a large enough constant to $v_n$, we can assume that the function $w_n:=-\log- v_n$ is qpsh and in $W^*$ : it satisfies $T_n:=dd^c w_n+ \omega \geq 0$ and $dw_n\wedge d^c w_n\leq T_n$ (see \cite{DS4}[Proposition 4.7] or \cite{moi2}). Since $w_n \geq v_n$, we have that $T_n$ is more H-regular than  $(1/d)^{n} (f^*)^{n}(\omega)$.
 Now for $A >0$, we consider the smooth function: 
$$\chi_A := h\left(\frac{w_n}{A}\right).$$
Then (up to substracting a large enough constant to $v_n$), there exists a constant $C$ that does not depend on $A$ such that:
\begin{itemize}
\item $d \chi_A \wedge d^c \chi_A = \frac{|h'(w_n/A)|^2}{A^2} d w_n \wedge d^c w_n \leq \frac{C}{A^2} T_n$.
\item $0 \leq  \frac{C}{A}(T_n+\omega) \pm dd^c\chi_A $.
\item $\chi_A =0$ in a neighborhood of $I(f^{n})$
\item $\lim_{A\to \infty} 1-\chi_A \searrow  1_{I(f^{n})}$ where $1_{I(f^{n})}$ is the indicator function of $I(f^{n})$.
\end{itemize}

\section{Construction of the measure}\label{Construction of the measure}
Let $\varphi \in C^1$. We define the quantity $c_n$ and the function $\varphi_n$ by:
$$\begin{cases} 
c_n:=\int_{\P^k} f^*(\varphi_{n-1}) T^+ _s\wedge \omega^{k-s}\\
\varphi_n:= f^*(\varphi_{n-1})-c_n
\end{cases}$$
with $c_0:=\int_{\P^k} \varphi T_s^+ \wedge \omega^{k-s}$ and $\varphi_0:= \varphi-c_0$. Observe that those quantities are well defined since $f^*(\varphi_{n-1})$ is smooth outside $\cap_{k\leq n} f^{-k}(I)$ which is an analytic set that has no mass for $T_s^+$. By construction, we have:
$$(f^n)^* \varphi = \sum_{i=0}^{n} c_i+ \varphi_n \ \mathrm{and} \  \int_{\P^k} \varphi \circ f^n T^+_s \wedge \omega^s =\sum_{i=0}^{n} c_i. $$
Hence, the following proposition shows in particular that $T_s^+ \wedge d^{-sn}f^n_*\omega^{k-s}$ converges to $\mu$ in the sense of currents with exponential estimates for $C^1$-observables.
\begin{proposition}\label{rate_measure}
There exists a constant $C$, independent of $n$ and $\varphi$ such that 
\begin{equation}\label{estime1}
 |c_n| \leq C \sqrt{\delta}^{-n} \|\varphi\|_{C^1}. 
\end{equation}
Consequently, $ |\mu(\varphi_n)|=|\mu(\varphi) -\sum_0^n c_i|\leq C \sqrt{\delta}^{-n} \|\varphi\|_{C^1}$.

Assume furthermore than $\varphi \in C^2$ and that $f$ is a birational regular map. Then
\begin{equation}\label{estime2}
 |c_n| \leq C \delta^{-n} \|\varphi\|_{C^2}. 
\end{equation}
Consequently, $  |\mu(\varphi_n)|=|\mu(\varphi) -\sum_0^n c_i|\leq C \delta^{-n} \|\varphi\|_{C^2}$.
\end{proposition}
\noindent \emph{Proof.} 
We explain in details the approximation's arguments as we will need them in several places. \\

\noindent \emph{Step 1.} The case $n=0$ is clear. We claim that: 
\begin{equation}\label{step1}
c_{n+1}= \int_{\P^k} \varphi_{n}T_s^+\wedge d_s^{-1}f_*(\omega^{k-s}).
\end{equation}
Assume first that $\varphi_{n}$ is continuous.  
We have that:
 $$c_{n+1}=\int_{\P^k} f^*(\varphi_{n}) T^+ _s\wedge \omega^{k-s} = \int_{\P^k} f^*(\varphi_{n})\frac{1}{d_s}f^*(T^+_s)\wedge \omega^{k-s}.$$
  Let $T_{1,\varepsilon}^+$ be a smooth approximation of $T^+_1$ (in the sense that it H-converges to $T_1^+$), then the quasi-potential $d^{-1} f^*(T_{1,\varepsilon}^+)$ also H-converges to $d^{-1} f^*(T_1^+)$. Then, $1/d_s f^*((T_{1,\varepsilon}^+)^s)$ converges to $1/d_sf^*(T^+_s)$ and thus $1/d_s f^*((T_{1,\varepsilon}^+)^s) \wedge \omega^{k-s}$ converges to  $ T^+ _s\wedge \omega^{k-s}$ in the sense of currents and $T^+_s\wedge \omega^{k-s}$ has no mass on $I(f)$ as $T_s^+$ does not. Similarly, as $T^+_s$ and $d_s^{-1}f_*(\omega)^{k-s}$ are wedgeable ($d_s^{-1}f_*(\omega)^{k-s}$ is more H-regular than $T^-_{k-s}$), $(T_{1,\varepsilon}^+)^s \wedge d_s^{-1}f_*(\omega)^{k-s}$ converges to  $T_{s}^+ \wedge 1/d_sf_*(\omega)^{k-s}$ in the Hartogs' sense (hence in the sense of currents).   
  In particular, since $f^*(\varphi_{n})\in L^1(T^+_s\wedge \omega^{k-s})$ and $\varphi_{n}$ is continuous, $c_{n+1}$ is the limit of 
\begin{align*}
\int_{\P^k\backslash I(f)} f^*(\varphi_{n}) 1/d_s f^*((T_{1,\varepsilon}^+)^s) \wedge \omega^{k-s}&= \int_{\P^k\backslash I(f)}  1/d_s f^*(\varphi_{n} (T_{1,\varepsilon}^+)^s) \wedge \omega^{k-s}\\
                                                                                       &=\int_{\P^k}  1/d_s f^*(\varphi_{n} (T_{1,\varepsilon}^+)^s) \wedge \omega^{k-s}\\
                  &=\int_{\P^k} \varphi_{n} (T_{1,\varepsilon}^+)^s \wedge \frac{1}{d_s}f_*(\omega)^{k-s}.
\end{align*} 
 Hence (\ref{step1}) stands for $\varphi_{n}$ continuous. An approximation of $\varphi_{n}$ gives the result by dominated convergence.\\
 
\noindent \emph{Step 2.} Let $T_{1,\varepsilon}^+$ be as in the previous paragraph. Let $S_\varepsilon$ be a smooth $H$-approximation of $\delta^{-1} f_*(\omega)=\omega+dd^c u$ so that $u_\varepsilon$ is a quasi-potential of $S_\varepsilon$ that decreases to $u$. We claim that:
\begin{align}\label{step2}
c_{n+1}=  \lim_{\varepsilon \to 0} \lim_{\varepsilon'\to 0}\int_{\P^k}  \varphi_{n} dd^c u_\varepsilon \wedge (T_{1,\varepsilon'}^+)^s \wedge R_\varepsilon
\end{align}
where $R_\varepsilon=\sum_{i=0}^{k-s-1}  S_\varepsilon^i\wedge \omega^{k-s-1-i} $.
 
Indeed, we write that: 
 $$ \frac{1}{d_s} f_*(\omega^{k-s}) = \omega^{k-s} + dd^c (u \wedge R)$$ 
where $R=\sum_{i=0}^{k-s-1} \delta^{-i} f_*(\omega)^i\wedge \omega^{k-s-1-i} $ is a positive closed current (of mass $k-s$). In particular, $u \wedge R$ is a quasi-potential of  
 $ 1/d_s f_*(\omega^{k-s})$. Similarly, $u_\varepsilon R_\varepsilon $ is a quasi-potential of $S_\varepsilon^{k-s}$ and $R_\varepsilon$ is a positive closed current (of mass $k-s$).  We have that $S_\varepsilon^{k-s}$ converges to $ 1/d_s f_*(\omega^{k-s})$ in the Hartogs' sense. 
 Continuity of the wedge product in that settings implies that $(T_{1,\varepsilon'}^+)^s\wedge S_\varepsilon^{k-s} $ converges to $T_s^+\wedge d_s^{-1}f_*(\omega^{k-s})$ in the sense of currents (and $R_\varepsilon \to R$ in the Hartogs' sense). In particular, since $\varphi_n \in L^\infty$ is smooth outside $I(f^n)$ and $I(f^n)$ has no mass for $T^+_s \wedge d_s^{-1} f_*(\omega_{k-s})$ (indeed, $T^+_s \wedge T^-_{k-s}$ gives no mass to $I(f^n)$ and $d_s^{-1} f_*(\omega_{k-s})$ is more H-regular than  $T^-_{k-s}$ hence neither does $T^+_s \wedge d_s^{-1} f_*(\omega_{k-s})$) :
\begin{align*}
c_{n+1} = \lim_{\varepsilon \to 0}\lim_{\varepsilon' \to 0} \int_{\P^k} \varphi_{n} (T_{1,\varepsilon'}^+)^s\wedge (\omega^{k-s} + dd^c (u_\varepsilon \wedge R_\varepsilon)) \\
 =  \lim_{\varepsilon \to 0} \lim_{\varepsilon' \to 0}\int_{\P^k} \varphi_{n} (T_{1,\varepsilon'}^+)^s\wedge \omega^{k-s}+
  \lim_{\varepsilon \to 0} \lim_{\varepsilon' \to 0}\int_{\P^k} \varphi_{n} (T_{1,\varepsilon'}^+)^s\wedge dd^c (u_\varepsilon \wedge R_\varepsilon).
\end{align*}   
The first term converges to $\int_{\P^k} \varphi_{n} T_s ^+\wedge \omega^{k-s} =0$ by construction. Now, as $ dd^c (u_\varepsilon \wedge R_\varepsilon)= dd^c (u_\varepsilon) \wedge R_\varepsilon$, step 2 follows. \\

\noindent \emph{Step 3.} We prove (\ref{estime1}). Observe that $ \varphi_{n} $ is in the Sobolev space $H^1$ since $d\varphi_n\wedge d^c \varphi_n = (f^{n})^* (d \varphi \wedge d^c \varphi) \leq \| \varphi\|^2_{C^1}(f^{n})^*( \omega)$ (this even implies that $\varphi_n \in W^*$ by \cite{DS4}). By Stokes, we can then write:
$$\int_{\P^k}  \varphi_{n} dd^c u_\varepsilon \wedge (T_{1,\varepsilon'}^+)^s \wedge R_\varepsilon =\int_{\P^k}  -d\varphi_{n} \wedge d^c u_\varepsilon \wedge (T_{1,\varepsilon'}^+)^s \wedge R_\varepsilon.$$
We write (the function $\chi_A$ was defined at the end of the previous section):
\begin{align*}
&\int_{\P^k}  -d\varphi_{n} \wedge d^c u_\varepsilon \wedge (T_{1,\varepsilon'}^+)^s \wedge R_\varepsilon= \\
&\int_{\P^k}  -\chi_Ad\varphi_{n} \wedge d^c u_\varepsilon \wedge (T_{1,\varepsilon'}^+)^s \wedge R_\varepsilon+ \int_{\P^k}  (\chi_A-1)d\varphi_{n} \wedge d^c u_\varepsilon \wedge (T_{1,\varepsilon'}^+)^s \wedge R_\varepsilon  \\
& = I + II
\end{align*}
with obvious notations. The first term converges, when $\varepsilon'\to 0$ to:
$$ \int_{\P^k}  -\chi_Ad\varphi_{n} \wedge d^c u_\varepsilon \wedge T_s^+ \wedge R_\varepsilon. $$
 By Cauchy-Schwarz inequality, we have that: 
\begin{align*}
&\Big|\int_{\P^k}  -\chi_Ad\varphi_{n} \wedge d^c u_\varepsilon \wedge T_s^+ \wedge R_\varepsilon\Big| \leq \\
&\left(\int_{\P^k} \chi_A d\varphi_{n} \wedge d^c \varphi_{n} \wedge T_s^+  \wedge R_\varepsilon \right)^{\frac{1}{2}}\left(\int_{\P^k} \chi_A du_\varepsilon \wedge d^c u_\varepsilon \wedge T_s^+\wedge R_\varepsilon \right)^{\frac{1}{2}}.
\end{align*}
Using $d\varphi_n\wedge d^c \varphi_n  \leq \| \varphi\|^2_{C^1}(f^{n})^*( \omega)$ and since  $f^{n}$ is smooth on the support of $\chi_A$, we have:  
\begin{align*}
 \chi_A d\varphi_{n} \wedge d^c \varphi_{n} \wedge T_s^+ &\leq \chi_A \| \varphi\|^2_{C^1} (f^{n})^*(\omega) \wedge \frac{1}{d_s^{n}}(f^{n})^*(T_s^+)\\
                                        &\leq \chi_A \| \varphi\|^2_{C^1} \frac{1}{d_s^{n}}(f^{n})^*(\omega \wedge T_s^+). 
\end{align*}
Then, we have that $\left(\int_{\P^k} \chi_A d\varphi_{n} \wedge d^c \varphi_{n} \wedge T_s^+  \wedge R_\varepsilon \right)^{\frac{1}{2}}$ is bounded by:
\begin{align*}
\left(\int_{\P^k} \chi_A  \| \varphi\|^2_{C^1} \frac{1}{d_s^{n}}(f^{n})^*(\omega \wedge T_s^+)  \wedge R_\varepsilon \right)^{\frac{1}{2}} &\leq \\
\| \varphi\|_{C^1} \left(\int_{\P^k}  \frac{d_{s+1}^n}{d_s^{n}d_{s+1}^n}(f^{n})^*(\omega \wedge T_s^+)  \wedge R_\varepsilon \right)^{\frac{1}{2}}&=\left(\frac{d_{s+1}}{d_s}\right)^{\frac{n}{2}}(k-s) \|\varphi\|_{C^1}\\
                     &\leq \sqrt{\delta}^{-n}  (k-s) \|\varphi\|_{C^1},
\end{align*}
as the last integral can be computed in cohomology ($k-s$ is the mass of $R_\varepsilon$ and $d_{s+1}= \delta^{k-s-1}$). The integral $ \int_{\P^k} \chi_A du_\varepsilon \wedge d^c u_\varepsilon \wedge T_s^+\wedge R_\varepsilon$ can be bounded by $\int_{\P^k} du_\varepsilon \wedge d^c u_\varepsilon \wedge T_s^+\wedge R_\varepsilon$ which by Cauchy-Schwarz inequality is just:
 $$\int_{\P^k} -u_\varepsilon \wedge dd^c u_\varepsilon \wedge T_s^+\wedge R_\varepsilon .$$
Since the quasi-potential $u$ is integrable with respect to $T^+_s \wedge 1/d_s f^{n}_*(\omega^{k-s})$ we have by Hartogs' convergence that that quantity converges (when $\varepsilon \to 0$) to the finite value $\int u dd^c u T^+_s \wedge R$  hence it is bounded for $\varepsilon$ small enough. In particular, $I$ is bounded by $C\|\varphi\|_{C^1} \sqrt{\delta}^{-n}$, which is exactly what we want.

We show now that $II$ can be taken arbitrarily small for $A$ large enough independent of $\varepsilon$, $\varepsilon'$. By Stokes:
\begin{align*}
II&=\int_{\P^k}  (1-\chi_A)\varphi_{n} \wedge dd^c u_\varepsilon \wedge (T_{1,\varepsilon'}^+)^s \wedge R_\varepsilon - \\
  &\  \int_{\P^k}   \varphi_{n} d\chi_A \wedge d^c u_\varepsilon \wedge (T_{1,\varepsilon'}^+)^s \wedge R_\varepsilon \\
  &=II'+II'',   
\end{align*}
again with obvious notations. Then:
\begin{align*}
II'=\int_{\P^k}  (1-\chi_A)\varphi_{n}  (T_{1,\varepsilon'}^+)^s \wedge S_\varepsilon\wedge R_\varepsilon- 
    \int_{\P^k}  (1-\chi_A)\varphi_{n}  (T_{1,\varepsilon'}^+)^s \wedge \omega \wedge R_\varepsilon.
\end{align*}
Since $\| \varphi_{n}\|_\infty$ is finite (say $\leq C_2$), we have the bound:
\begin{align*}
|II'|\leq C_2( \int_{\P^k}  (1-\chi_A) (T_{1,\varepsilon'}^+)^s \wedge S_\varepsilon\wedge R_\varepsilon+ 
    \int_{\P^k}  (1-\chi_A) (T_{1,\varepsilon'}^+)^s \wedge \omega \wedge R_\varepsilon).
\end{align*} 
 The integral $\int_{\P^k}  (1-\chi_A) (T_{1,\varepsilon'}^+)^s \wedge S_\varepsilon\wedge R_\varepsilon$ converges when $\varepsilon, \varepsilon' \to 0$ to $\int_{\P^k}  (1-\chi_A) T_s^+ \wedge \delta^{-1}f_*(\omega)\wedge  R$ which can be taken arbitrarily small for $A$ large enough (recall that $T_s
^+ \wedge \delta^{-1}f_*(\omega)\wedge R$ gives no mass to $I(f^{n})$). For the same reasons, the second integral can be taken arbitrarily small for $A$ large enough.
 
 Now the term $II''$ can be bounded by Cauchy-Schwarz inequality:
\begin{align*}
|II''| &\leq \left(\int_{\P^k}  d \chi_A \wedge d^c \chi_A \wedge (T_{1,\varepsilon'}^+)^s \wedge R_\varepsilon \right)^{\frac{1}{2}} 
           \left(\int_{\P^k}  d u_\varepsilon \wedge d^c u_\varepsilon \wedge (T_{1,\varepsilon'}^+)^s \wedge R_\varepsilon \right)^{\frac{1}{2}} \\
       &\leq  \left(\int_{\P^k}  \frac{1}{A^2} T_n \wedge (T_{1,\varepsilon'}^+)^s \wedge R_\varepsilon \right)^{\frac{1}{2}}  \left(\int_{\P^k}  d u_\varepsilon \wedge d^c u_\varepsilon \wedge (T_{1,\varepsilon'}^+)^s \wedge R_\varepsilon \right)^{\frac{1}{2}} 
\end{align*}
 where $T_n$ was defined in the paragraph defining $\chi_A$. In particular, the first integral can be computed in cohomology and is equal to $A^{-2}(k-s)$ hence it tends to $0$ (uniformly in  $\varepsilon, \ \varepsilon'$). The second was bounded in the previous paragraph.

 Summing up, we get that there exists a constant $C$ that does not depend on $n$ or $\varphi$ such that $|c_n| \leq C \sqrt{\delta}^{-n} \|\varphi\|_{C^1}$. This gives (\ref{estime1}) and the speed of convergence $ |\mu(\varphi) -\sum_0^n c_i|\leq C \sqrt{\delta}^{-n} \|\varphi\|_{C^1}$ follows by classical series' arguments ($\mu(\varphi)=\sum c_i$ hence $\mu(\varphi)-\sum_{i=0}^n c_i= \sum_{i\geq n+1} c_i$). \\

 \noindent \emph{Step 3'.} Assume now that $f$ is a regular birational maps of $\P^k$ and that $\varphi \in C^2$. We have that $\|\varphi\|_{C^2} \omega \pm dd^c\varphi \geq 0$, hence  $\|\varphi\|_{C^2} (f^{n})^*(\omega) \pm dd^c\varphi_{n} \geq 0$. It implies that $\varphi_n$ is dsh (difference of qpsh functions). Then Stokes formula implies: 
 $$\int_{\P^k}  \varphi_{n} dd^c u_\varepsilon \wedge (T_{1,\varepsilon'}^+)^s \wedge R_\varepsilon =\int_{\P^k}   u_\varepsilon  dd^c\varphi_{n} \wedge (T_{1,\varepsilon'}^+)^s \wedge R_\varepsilon.$$
 As in Step 3, we write it as: 
\begin{align*}
 \int_{\P^k}  \chi_A  u_\varepsilon  dd^c\varphi_{n} \wedge (T_{1,\varepsilon'}^+)^s \wedge R_\varepsilon  +   \int_{\P^k}  (1-\chi_A)  u_\varepsilon  dd^c\varphi_{n} \wedge (T_{1,\varepsilon'}^+)^s \wedge R_\varepsilon.
\end{align*}
When $\varepsilon' \to 0$, the first integral converges to (every term is smooth):
 $$ \int_{\P^k}  \chi_A  u_\varepsilon  dd^c\varphi_{n} \wedge T_{s}^+ \wedge R_\varepsilon.$$
 By invariance, $T_{s}^+=d_s^{-n}(f^{n})^*(T^+_s)$ and since $f^{n}$ is smooth in the support of $\chi_A$, we have that:
\begin{align*}
 0\leq  -\chi_A  u_\varepsilon  \|\varphi\|_{C^2} d_s^{-n}(f^{n})^*(\omega \wedge T_{s}^+) \wedge R_\varepsilon   \pm  \chi_A  u_\varepsilon  dd^c\varphi_{n} \wedge T_{s}^+ \wedge R_\varepsilon 
\end{align*}
 (recall that $u_\varepsilon <0$). Hence, we have the bound:
 $$  \left|\int_{\P^k}  \chi_A  u_\varepsilon  dd^c\varphi_{n} \wedge T_{s}^+ \wedge R_\varepsilon \right| \leq -\int_{\P^k}         \|\varphi\|_{C^2} d_s^{-n}(f^{n})^*(\omega \wedge T_{s}^+) \wedge u_\varepsilon R_\varepsilon .$$
Then the previous quantity can be written as: 
$$      - \|\varphi\|_{C^2} \frac{1}{\delta^{n}} \U_\varepsilon \left( d_{s+1}^{-n}(f^{n})^*(\omega \wedge T_{s}^+) \right)$$
 where $\U_\varepsilon$ is the super-potential of $S_\varepsilon^{k-s}$ associated to the quasi-potential $u_\varepsilon R_\varepsilon$. Hence, by continuity for the Hartogs' convergence, it converges when $\varepsilon \to 0$ to: 
 $$      - \|\varphi\|_{C^2} \frac{1}{\delta^{n}} \U \left(d_{s+1}^{-n}(f^{n})^*(\omega \wedge T_{s}^+) \right)  $$
 where $\U$ is the super-potential of $d_s^{-1} f_*(\omega^{k-s})$ associated to the quasi-potential $uR$. Recall that by \cite{DS10}, $\U$ is uniformly bounded at any positive closed current of mass $1$ and support in $\supp(T_s^+)$ which is the case of $d_{s+1}^{-n}(f^{n})^*(\omega \wedge T_{s}^+)$. This bounds the previous term by $C \|\varphi\|_{C^2} \delta^{-n}$ where $C$ does not depend on $n$ nor $\varphi$.

Now we show that the second integral $\int_{\P^k}  (1-\chi_A)  u_\varepsilon  dd^c\varphi_{n} \wedge (T_{1,\varepsilon'}^+)^s \wedge R_\varepsilon$ can be taken arbitrarily small for $A$ large enough. We could use simpler arguments (using the fact that $u$ is bounded near $I(f^{n})$) but we will need the following arguments later on hence observe that we will not use that $f$ is regular here.
 As in Step 3, we use Stokes  formula to write it as:
$$\int_{\P^k}   u_\varepsilon d\chi_A\wedge  d^c\varphi_{n} \wedge (T_{1,\varepsilon'}^+)^s \wedge R_\varepsilon + \int_{\P^k}  (\chi_A-1)  du_\varepsilon \wedge  d^c\varphi_{n} \wedge (T_{1,\varepsilon'}^+)^s \wedge R_\varepsilon.$$
Consider the first integral, by Cauchy-Schwarz inequality:
\begin{align*}
&\left| \int_{\P^k}   u_\varepsilon d\chi_A\wedge  d^c\varphi_{n} \wedge (T_{1,\varepsilon'}^+)^s \wedge R_\varepsilon \right|^2 \leq  \left( \int_{\P^k}   -u_\varepsilon d\chi_A\wedge  d^c\chi_A \wedge (T_{1,\varepsilon'}^+)^s \wedge R_\varepsilon \right) \\
& \left( \int_{\supp \chi_A}   -u_\varepsilon d\varphi_{n}\wedge  d^c\varphi_{n} \wedge (T_{1,\varepsilon'}^+)^s \wedge R_\varepsilon \right).
\end{align*}
The first integral of the product can be bounded by:
$$\int_{\P^k}   \frac{-u_\varepsilon}{A^2} T_n \wedge (T_{1,\varepsilon'}^+)^s \wedge R_\varepsilon  $$
which converges to:
$$\int_{\P^k}   \frac{-u}{A^2} T_n \wedge T_{s}^+ \wedge R .$$
when $\varepsilon, \ \varepsilon' \to 0$ which is finite since $u\in L^1(T^+_1\wedge T^+_s \wedge (T_1^-)^{s-1})$ (
$T_n$ is more H-regular than $(1/d)^{n} (f^*)^{n}(\omega)$ hence than $T^+_1$). Finally, it can be taken arbitrarily small for $A$ large enough. Similarly, since $d\varphi_{n}\wedge  d^c\varphi_{n} \leq  \| \varphi\|_{C^1}d^{-n}(f^{n})^*(\omega)$ and using the fact that $u\in L^1(T^+_1\wedge T^+_s \wedge (T_1^-)^{s-1})$ 
we have that the second integral of the product is bounded. In particular, the integral:
$$\int_{\P^k}   u_\varepsilon d\chi_A\wedge  d^c\varphi_{n} \wedge (T_{1,\varepsilon'}^+)^s \wedge R_\varepsilon $$
can be taken arbitrarily small for $A$ large enough.

To conclude, we have to show that  $\int_{\P^k}  (\chi_A-1)  du_\varepsilon \wedge  d^c\varphi_{n} \wedge (T_{1,\varepsilon'}^+)^s \wedge R_\varepsilon $ can also be taken arbitrarily small. We recognize that it is the same integral than the integral $II$ in Step 3 (up to exchanging $d$ and $d^c$). Hence it goes to $0$ which proves (\ref{estime2}) (again $ |\mu(\varphi) -\sum_0^n c_i|\leq C \delta^{-n} \|\varphi\|_{C^2}$ follows directly). \hfill $\Box$ \hfill
 
\section{Decay of correlations}\label{Decayofcorr}

We now prove:
\begin{proposition}
Let $\varphi$ and $\psi$ in $C^2$, $0\leq \alpha \leq 2$, then there exists a constant $C$ that does not depend on $\varphi$ and $\psi$ such that:
\begin{equation}\label{eq1}
\forall n, \ m, \  |\mu(\varphi \circ f^n . \psi \circ f^{-m}) - \mu(\varphi).\mu( \psi)|\leq C \|\varphi \|_{C^2}\|\psi \|_{C^2} (\sqrt{\delta}^{-n}+\sqrt{d}^{-m})  . 
\end{equation}
Furtheremore, if $f$ is a regular birational map then we have the estimate:
\begin{equation}\label{eq2}
\forall n, \ m, \  |\mu(\varphi \circ f^n . \psi \circ f^{-m}) - \mu(\varphi).\mu( \psi)|\leq C \|\varphi \|_{C^2}\|\psi \|_{C^2} (\delta^{-n}+d^{-m}). 
\end{equation}
\end{proposition}
\noindent \emph{Proof.} We consider the general case first and we will say after what are the easy modifications when $f$ is regular.  

We apply the results of the previous section and decompose $\varphi \circ f^n$ into $c_0 + \dots +c_n + \varphi_n$. Replacing $f$ by $f^{-1}$, we can write $\psi \circ f^{-m}$ as $c'_0 + \dots +c'_m + \psi_m$ where $\int \psi_m \omega^{s} \wedge T^-_{k-s}=0$. Then, we have by  Proposition \ref{rate_measure} that $c'_i \leq C \sqrt{d}^{-j} \|\psi\|_{C^1}$ and $ |\mu(\psi_m)|=|\mu(\psi) -\sum_0^m c'_i|\leq C \sqrt{d}^{-m} \|\psi\|_{C^1}$.

In particular, we have that ($\mu$ is invariant):
$$\mu(\varphi)=\sum_0^n c_i + \mu(\varphi_n), \quad    \mu(\psi)=\sum_0^n c'_i + \mu(\psi_m).$$
Hence:
 \begin{align*} \mu(\varphi \circ f^n . \psi \circ f^{-m})&=\mu(  (\sum_0^n c_i +\varphi_n)( \sum_0^n c'_i +\psi_m  ))\\
                                                          &= \mu((\mu(\varphi)-\mu(\varphi_n)+\varphi_n).(\mu(\psi)-\mu(\psi_m)+\psi_m)).
  \end{align*}
After simplifications, we deduce:
\begin{align*}
\mu(\varphi \circ f^n . \psi \circ f^{-m}) - \mu(\varphi).\mu( \psi)=\mu(\varphi_n \psi_m)-\mu(\varphi_n)\mu(\psi_m).
  \end{align*}
  By Proposition \ref{rate_measure}, we have that $|\mu(\varphi_n)\mu(\psi_m)|\leq C^2 \|\varphi\|_{C^1}\|\psi\|_{C^1} \sqrt{d}^{-m}\sqrt{\delta}^{-n}$ (which is indeed $\leq C \|\varphi \|_{C^2}\|\psi \|_{C^2} (\sqrt{\delta}^{-n}+\sqrt{d}^{-m})$ where $C$ does not depend on $n$, $m$ nor the test functions). So all there is left to do is bounding  $\mu(\varphi_n \psi_m)$.
 We write:
$$\mu(\varphi_n . \psi_m)=\int_{(\P^k)^2} \varphi_n(x) \psi_m(y) T^+_s(x)\wedge T^-_{k-s}(y)\wedge [\Delta],  $$
  where $x$ denotes the coordinate on the first $\P^k$, $y$ the coordinate on the second $\P^k$ and  $[\Delta]$ is the current of integration on the diagonal $\Delta$ of $(\P^k)^2$. We write $[\Delta]= \sum_{i=0}^k \omega^i(x)\wedge \omega^{k-i}(y)+ dd^c V$ where $V$ is the negative quasi-potential of $\Delta$ given in \cite{DS6}[Theorem 2.3.1]. Bidegree's arguments imply that:
\begin{align*}
\mu(\varphi_n . \psi_m)&=\int_{(\P^k)^2} \varphi_n(x) \psi_m(y) T^+_s(x)\wedge T^-_{k-s}(y)\wedge \omega^{k-s}(x) \wedge \omega^s(y) \\
 & + \int_{(\P^k)^2} \varphi_n(x) \psi_m(y) T^+_s(x)\wedge T^-_{k-s}(y)\wedge dd^c V.
 \end{align*}     
By Fubini, the first integral is $0$ (recall that $\int \varphi_n T_s^+ \wedge\omega^{k-s}=0$ by definition). Let $\chi_A$ be the cut-off function defined earlier and let $\xi_A$ be the similar cut-off function associated to $I(f^{-n})$. Then we have that:
\begin{align*}
\mu(\varphi_n . \psi_m)= \lim_{A\to \infty} \int_{(\P^k)^2} \chi_A(x)\varphi_n(x) \xi_A(y)\psi_m(y) T^+_s(x)\wedge T^-_{k-s}(y)\wedge dd^c V.
 \end{align*}   
 Now, let $T^+_{s,\varepsilon}$ and $ T^-_{k-s, \varepsilon}$ be Hartogs' approximations of $T^+_s$ and $T^-_{k-s}$, then:
 \begin{align*}
\mu(\varphi_n . \psi_m)= \lim_{A\to \infty} \lim_{\varepsilon \to 0} \int_{(\P^k)^2} &\chi_A(x)\varphi_n(x) \xi_A(y)\psi_m(y) \frac{1}{d_s^{n}}(f^{n})^*T^+_{s,\varepsilon}(x) \\
&\wedge \frac{1}{d_s^{m}}(f^{m})_*T^-_{k-s,\varepsilon}(y)\wedge dd^c V.
 \end{align*}  
 We apply Stokes formula twice, this gives sixteen terms (we left the writing of the full expression to the reader). Let us show that every one of these terms with a $d \chi_A$, $d^c \chi_A$, or $dd^c \chi_A$ can be taken arbitrarily small for $A$ large enough independently of $\varepsilon$ (same thing for $\xi_A$).   First, we consider the terms: 
 \begin{align}\label{exp5} 
 \int_{(\P^k)^2} &\xi_A(y)\psi_m(y) (d\chi_A(x)\wedge d^c\varphi_n(x) + d\varphi_n(x)\wedge d^c\chi_A(x)) \\
  &\wedge  \frac{1}{d_s^{n}}(f^{n})^*T^+_{s,\varepsilon}(x) \wedge \frac{1}{d_s^{m}}(f^{m})_*T^-_{k-s,\varepsilon}(y)\wedge  V  \notag .
  \end{align}
Using Cauchy-Schwarz inequality (see \cite{dem}[Chapter III] for bases on (strongly) positive currents), we have that its square is less than:
 \begin{align}\label{expression1} 
 &C_1 \int_{(\P^k)^2} - d\varphi_n(x)\wedge d^c\varphi_n(x) \wedge  \frac{1}{d_s^{n}}(f^{n})^*T^+_{s,\varepsilon}(x) \wedge \frac{1}{d_s^{m}}(f^{m})_*T^-_{k-s,\varepsilon}(y)\wedge  V \\
  & \times \int_{(\P^k)^2} -d\chi_A(x)\wedge d^c\chi_A(x) \wedge  \frac{1}{d_s^{n}}(f^{n})^*T^+_{s,\varepsilon}(x) \wedge \frac{1}{d_s^{m}}(f^{m})_*T^-_{k-s,\varepsilon}(y)\wedge  V \notag
  \end{align}  
 where $C_1=\|\psi_m\|^2_\infty$. 
We show that the first term of the product is uniformly bounded for $\varepsilon$ small enough. As in the previous section, we use that  $ d\varphi_{n} \wedge d^c \varphi_n \leq \|\varphi\|^2_{C^1} (f^{n})^*\omega$. Observe that the integral can be taken outside $I(f^{n})$, so the term is bounded by:
$$ - \int_{(\P^k)^2} \|\varphi\|^2_{C^1} \frac{1}{d_s^{n}}(f^{n})^*(\omega \wedge T^+_{s,\varepsilon}(x)) \wedge \frac{1}{d_s^{m}}(f^{m})_*T^-_{k-s,\varepsilon}(y)\wedge  V .$$
Now, if $p_1$ and $p_2$ denote the projections from $(\P^k)^2$ onto its factors and $R$ is a positive closed current on $\P^k$, we have that $(p_1)_*(p_2^*(R)\wedge V)$ is the Green quasi-potential of $R$ (see \cite{DS6}[Theorem 2.3.1]) that we denote by $U_{R}$. Here, for $R=\frac{1}{d_s^{m}}(f^{m})_*T^-_{k-s,\varepsilon}$, we have that the term is bounded by:
 $$ -\|\varphi\|^2_{C^1} \int_{\P^k}  \frac{1}{d_s^{n}}(f^{n})^*(\omega \wedge T^+_{s,\varepsilon}) \wedge U_{\frac{1}{d_s^{m}}(f^{m})_*T^-_{k-s,\varepsilon}} .$$  
 The following lemma will give us the needed estimate.
\begin{lemme}\label{keyCS}
There exists a constant $K$ independent of $n$ such that for $\varepsilon$ small enough:
\begin{align}\label{eq3}  \left| \int_{\P^k} \frac{1}{d_{s+1}^{n}}(f^{n})^*(\omega\wedge T^+_{s,\varepsilon})\wedge U_{\frac{1}{d_s^{m}}(f^{m})_*T^-_{k-s,\varepsilon} }     \right| \leq K  \sqrt{\delta}^{n} ,
\end{align} 
and:
\begin{align}\label{eq4}  \left| \int_{\P^k} \frac{1}{d_{s-1}^{n}}(f^{n})_*(\omega\wedge T^-_{k-s,\varepsilon})\wedge U_{\frac{1}{d_s^{n}} (f^{n})^* T^+_{s,\varepsilon} }     \right| \leq K  \sqrt{d}^n. 
\end{align} 
\end{lemme}
\noindent \emph{Proof of the lemma.} In terms of super-potentials, the integral in (\ref{eq3}) can be written as:
$$  \U_{\frac{1}{d_s^{m}}(f^{m})_*T^-_{k-s,\varepsilon} }\left( \frac{1}{d_{s+1}^{n}}(f^{n})^*(\omega\wedge T^+_{s,\varepsilon})\right) .$$
The current $T^+_{s,\varepsilon}$ converges to $T^+_{s}$ in the Hartogs' sense hence $\omega\wedge T^+_{s,\varepsilon}$ converges to $\omega\wedge T^+_{s} $ in the Hartogs' sense. The current $\omega\wedge T^+_{s}$ is $(f^{n})^*$-admissible since its super-potential is finite at $T^-_{k-s}= d_s^{-n}(f^{n})_*(T^-_{k-s})$. Since,  $T^-_{k-s,\varepsilon} \to T^-_{k-s}$ in the Hartogs' sense and $T^-_{k-s}$ is $f^{n}_*$-admissible, we have that when $\varepsilon \to 0$: 
$$\frac{1}{d_s^{m}}(f^{m})_*T^-_{k-s,\varepsilon}\to \frac{1}{d_s^{m}}(f^{m})_*T^-_{k-s}= T^-_{k-s}$$
in the Hartogs' sense.
Now, continuity of pull-back and evaluation for the Hartogs' convergence implies that: 
 $$ \U_{\frac{1}{d_s^{m}}(f^{m})_*T^-_{k-s,\varepsilon} } (\frac{1}{d_{s+1}^{n}}(f^{n})^*(\omega\wedge T^+_{s,\varepsilon}) ) \to   \U_{T^-_{k-s} }(  \frac{1}{d_{s+1}^{n}}(f^{n})^*(\omega\wedge T^+_{s})).  $$
Since $T^-_1 $ and $T^+_s$ are wedgeable and $T^-_1 \wedge T^+_s$ is a $(f^{n})^*$-invariant current  (\cite{DV1}[Theorem 3.3.8]), we have that 
 $\U_{T^-_{k-s} }(  \frac{1}{d_{s+1}^{n}}(f^{n})^*(T^-_1\wedge T^+_{s}))$ is well defined and equal to $\U_{T^-_{k-s} }(T^-_1\wedge T^+_{s})$ (hence it is bounded). By difference, we have to control (we extend the definition of super-potentials by linearity to linear combinations of positive closed currents):
 $$ \U_{T^-_{k-s} }(  \frac{1}{d_{s+1}^{n}}(f^{n})^*( dd^c U_{T^-_1} \wedge T^+_{s})).$$
 We can choose a particular super-potential of $T^-_{k-s}$ (that will just change the result by an additive constant). So we choose the super-potential defined by the quasi-potential $U_{T^-_1}\wedge Q$ where $U_{T^-_1}$ is a quasi-potential of $T^-_1$ and  $Q$ is the current  $Q=\sum_{i=0}^{k-s-1} (T^{-}_1)^i\wedge \omega^{k-s-1-i}$ (recall that $T^{-}_{k-s}= (T^-_1)^{k-s}$ in the sense of super-potentials and use \cite{DV1}[Lemma A.2.1]). So we have to control (everything can be written in terms of $(1,1)$ currents and qpsh functions):
 $$  \langle U_{T^-_1}\wedge Q , \frac{1}{d_{s+1}^{n}}(f^{n})^*( dd^c U_{T^-_1} \wedge T^+_{s}) \rangle.$$
 Now, in order to do the following computations, one should, exactly as in the previous section (Step III'), H-regularize every term in order to deal with smooth forms and use the cut-off function $\chi_A$ to show that nothing happen on $I(f^{n})$. Though we will not do that in order to simplify the exposition, it is what we do implicitly, the arguments being exactly the same than in the previous section. Now using Stokes formula, we write the previous term as:
  $$  -\langle d U_{T^-_1}\wedge Q , \frac{1}{d_{s+1}^{n}}(f^{n})^*(d^c U_{T^-_1} \wedge T^+_{s}) \rangle.$$
 Using the invariance of $T^+_s$, it is:
   $$  -\langle d U_{T^-_1}\wedge Q , \frac{d_s^n}{d_{s+1}^{n}}d^c(f^{n})^*( U_{T^-_1}) \wedge T^+_{s} \rangle.$$
 Now, applying Cauchy-Schwartz inequality, it is bounded by:
 \begin{align*}
  & \frac{d_s^n}{d_{s+1}^{n}} \left(\int d U_{T^-_1}\wedge d^c U_{T^-_1} \wedge Q \wedge T^+_{s}\right)^{\frac{1}{2}} \times \\
    &    \left( \int d(f^{n})^*( U_{T^-_1}) \wedge d^c(f^{n})^*( U_{T^-_1})\wedge Q  \wedge T^+_{s}\right)^{\frac{1}{2}} .
 \end{align*}
  The first integral is bounded since the quasi-potential of $T^-_1$ is integrable with respect to $\mu=T^+_s\wedge T^-_{k-s}$ and $Q$ is a combination of terms more H-regular than $(T^-_1)^{s-1}$. The second integral can be written as:
  $$\int d(f^{n})^* U_{T^-_1} \wedge d^c(f^{n})^* U_{T^-_1}  \wedge T^+_{s}\wedge Q= \int \frac{1}{d_s^{n}}(f^{n})^*( dU_{T^-_1} \wedge d^cU_{T^-_1}  \wedge T^+_{s})\wedge Q,$$
this is where one should use the cut-off function $\chi_A$ to show that there is no mass on $I(f^{n})$ and then use that $d_s^{-n}(f^{n})^*(T^+_s)=T^+_s$. Pushing forward and using Stokes formula, we write it as:
$$  \int \frac{d_{s+1}^{n}}{d_s^{n}}( -U_{T^-_1} \wedge dd^cU_{T^-_1}  \wedge T^+_{s})\wedge  \frac{1}{d_{s+1}^{n}}(f^{n})_*Q. $$ 
 Now as $d_{s+1}=\delta^{k-s-1}$, we have that: 
 $$\frac{1}{d_{s+1}^{n}}(f^{n})_*Q=\sum_{i=0}^{k-s-1} (T^{-}_1)^i\wedge \frac{1}{\delta^{(k-s-1-i)n}}(f^{n})_*\omega^{k-s-1-i}$$  
Since:
 $$(T^{-}_1)^i\wedge \frac{1}{\delta^{(k-s-1-i)n}}(f^{n})_*\omega^{k-s-1-i}\to T^-_{k-s-1}=(T^-_1)^{k-s-1}$$
  in the Hartogs' sense and $U_{T^-_1}$ is integrable with respect to $T^+_s\wedge T^-_{k-s}$, each  quantity  $  \int  -U_{T^-_1} \wedge dd^cU_{T^-_1}  \wedge T^+_{s}\wedge  (T^{-}_1)^i\wedge \frac{1}{\delta^{(k-s-1-i)n}}(f^{n})_*\omega^{k-s-1-i}$ is uniformly bounded by a constant that does not depend on $n$. Hence: 
 $$  \int -U_{T^-_1} \wedge dd^cU_{T^-_1}  \wedge T^+_{s}\wedge  \frac{1}{d_{s+1}^{n}}(f^{n})_*Q $$  
 is  bounded by a constant that does not depend on $n$. Using that  $d_s^n= \delta^{n}d_{s+1}^{n}$
  gives (\ref{eq3}), the proof of (\ref{eq4}) is the same. \hfill $\Box$ \hfill  \\

In particular, the first term in the product of (\ref{expression1}) is bounded and using that $d_{s}=\delta d_{s+1}$ we have proved:
 \begin{align}\label{boundkey} 
 &\int_{(\P^k)^2} - d\varphi_n(x)\wedge d^c\varphi_n(x) \wedge  \frac{1}{d_s^{n}}(f^{n})^*T^+_{s,\varepsilon}(x) \wedge\frac{1}{d_s^{m}}(f^{m})_*T^-_{k-s,\varepsilon}(y)\wedge  V \notag\\
 &  \leq  \|\varphi\|^2_{C^1} \frac{K}{\sqrt{\delta}^n}                 
   \end{align}
(we will need the above precise bound later). We show now that the second term of the product in (\ref{expression1}), denoted by $I_A$, goes to $0$ when $A\to \infty$. Recall that the function $\chi_A$ satisfies  $d \chi_A \wedge d^c \chi_A \leq \frac{1}{A^2} T_n$ where $T_n= dd^cw_n+\omega$ (see the paragraph at the end of Section \ref{setting}). Hence $I_A$ is bounded by:
$$ \int_{(\P^k)^2}  -\frac{1}{A^2}T_n\wedge \frac{1}{d_s^{n}}(f^{n})^*(T^+_{s,\varepsilon}(x)) \wedge \frac{1}{d_s^{m}}(f^{m})_*T^-_{k-s,\varepsilon}(y)\wedge  V.$$
 As above, we write it as:
$$ \int_{\P^k} -\frac{1}{A^2}T_n\wedge \frac{1}{d_s^{n}}(f^{n})^*(T^+_{s,\varepsilon})\wedge U_{\frac{1}{d_s^{m}}(f^{m})_*T^-_{k-s,\varepsilon}}. $$
In terms of super-potentials, the previous quantity can be interpreted as:
$$  \frac{-1}{A^2} \U_{\frac{1}{d_s^{m}}(f^{m})_*T^-_{k-s,\varepsilon}}\left(T_n\wedge \frac{1}{d_s^{n}}(f^{n})^*(T^+_{s,\varepsilon}(x))\right),$$
where $ \U_{\frac{1}{d_s^{m}}(f^{m})_*T^-_{k-s,\varepsilon}}$ is the super-potential of $\frac{1}{d_s^{m}}(f^{m})_*T^-_{k-s,\varepsilon}$ associated to the Green quasi-potential $U_{\frac{1}{d_s^{m}}(f^{m})_*T^-_{k-s,\varepsilon}}$. Since $T_n$ is more H-regular than $T^+_1$  and  
since $T^+_1$ and $T^+_{s}$ are wedgeable we can use the same arguments than in Lemma \ref{keyCS} to show that the above quantity converges to:
$$\U_{T^-_{k-s}}(T_n\wedge T^+_s)>-\infty.$$ 
In particular,  dividing by $A^2$ gives that the term $I_A$ can be taken arbitrarily small for $A$ large enough independently of $\varepsilon$.
So, the  term (\ref{exp5}) can be taken arbitrarily small for $A$ large enough. \\   
  
 The same computations give that the terms in $d\chi_A \wedge d^c\psi_m + d\psi_m  \wedge d^c\chi_A $, $d\xi_A \wedge d^c\varphi_n+ d\varphi_n  \wedge d^c\xi_A$, $d\xi_A \wedge d^c\psi_m+ d\psi_m  \wedge d^c\xi_A$ can all be taken arbitrarily small for $A$ large enough. Now, we consider the term:
 \begin{align}\label{exp6} 
 \int_{(\P^k)^2} &\xi_A(y)\psi_m(y)\varphi_n(x) dd^c\chi_A(x)\wedge  \frac{1}{d_s^{n}}(f^{n})^*T^+_{s,\varepsilon}(x) \wedge \frac{1}{d_s^{m}}(f^{m})_*T^-_{k-s,\varepsilon}(y)\wedge  V   .
  \end{align}
 Recall that $0 \leq  \frac{C}{A}(T_n+\omega) \pm dd^c\chi_A $ where $C$ is another constant that does not depend on $A$. Hence:
$$0 \leq  \frac{C'}{A}(T_n+\omega) \pm  \xi_A(y)\psi_m(y)\varphi_n(x) dd^c\chi_A(x) dd^c\chi_A $$
 where $C'=\|\psi_m\|_\infty \|\varphi_n\|_\infty C$ is a constant that does not depend on $A$. Hence the term (\ref{exp6}) is bounded in absolute value by:
 $$ \frac{C'}{A} \int_{(\P^k)^2} -(T_n+\omega)\wedge \frac{1}{d_s^{n}}(f^{n})^*(T^+_{s,\varepsilon}(x)) \wedge \frac{1}{d_s^{m}}(f^{m})_*T^-_{k-s,\varepsilon}(y)\wedge  V.$$
 The same argument that above ($T_n+\omega$ is more H-regular than $T^+_1$) implies that this term can be taken arbitrarily small for $A$ large enough independently of $n$ and $m$. By symmetry, the term in $dd^c\xi_A$ can also be taken arbitrarily small for $A$ large. \\
 
  So all there is left to bound are the following terms:
  \begin{align*} 
 \int_{(\P^k)^2} &\chi_A(x)\xi_A(y) (d\varphi_n(x)\wedge d^c\psi_m(y)+ d\psi_m(y)\wedge d^c\varphi_n(x))\wedge  \frac{1}{d_s^{n}}(f^{n})^*T^+_{s,\varepsilon}(x)\\
 & \wedge \frac{1}{d_s^{m}}(f^{m})_*T^-_{k-s,\varepsilon}(y)\wedge  V,  \\
  \int_{(\P^k)^2} &\chi_A(x)\xi_A(y)\psi_m(y) dd^c\varphi_n(x)\wedge  \frac{1}{d_s^{n}}(f^{n})^*T^+_{s,\varepsilon}(x) \wedge \frac{1}{d_s^{m}}(f^{m})_*T^-_{k-s,\varepsilon}(y)\wedge  V , \\
    \int_{(\P^k)^2}& \chi_A(x)\xi_A(y)\varphi_n(y) dd^c\psi_m(y)\wedge  \frac{1}{d_s^{n}}(f^{n})^*T^+_{s,\varepsilon}(x) \wedge \frac{1}{d_s^{m}}(f^{m})_*T^-_{k-s,\varepsilon}(y)\wedge  V . \\
      \end{align*} 
 The first one of these term is bounded by Cauchy-Schwarz inequality by:
    \begin{align*} 
 \left( \int_{(\P^k)^2} d\varphi_n(x)\wedge d^c\varphi_n(x)\wedge  \frac{1}{d_s^{n}}(f^{n})^*T^+_{s,\varepsilon}(x) \wedge \frac{1}{d_s^{m}}(f^{m})_*T^-_{k-s,\varepsilon}(y)\wedge  V \right)^\frac{1}{2}\times\\
  \left( \int_{(\P^k)^2} d\psi_m(x)\wedge d^c\psi_m(x)\wedge  \frac{1}{d_s^{n}}(f^{n})^*T^+_{s,\varepsilon}(x) \wedge \frac{1}{d_s^{m}}(f^{m})_*T^-_{k-s,\varepsilon}(y)\wedge  V \right)^\frac{1}{2}
      \end{align*}  
We have already bounded the first term of the product in (\ref{boundkey}), we have a similar bound for the other term. Hence, the product is bounded by:
$$  K^2 \|\varphi\|_{C_1}\|\psi\|_{C_1}  \sqrt{d}^{-m} \sqrt{\delta}^{-n} \leq K^2\|\varphi\|_{C_2}\|\psi\|_{C_2}( \sqrt{\delta}^{-n} + \sqrt{d}^{-m}).$$ 
We now control the term in $dd^c\varphi_n(x)$. We use that $\|\varphi\|_{C^2} (f^{n})^*(\omega) \pm dd^c\varphi_{n} \geq 0$ (see Step 3' in previous section). It implies that: 
$$\|\psi\|_\infty \|\varphi\|_{C^2} (f^{n})^*(\omega) \pm \chi_A(x)\xi_A(y)\psi_m(y)dd^c\varphi_{n} \geq 0.$$
 That gives:
   \begin{align*} 
 \Big| \int_{(\P^k)^2} &\chi_A(x)\xi_A(y)\psi_m(y) dd^c\varphi_n(x)\wedge  \frac{1}{d_s^{n}}(f^{n})^*T^+_{s,\varepsilon}(x) \\
  &\wedge \frac{1}{d_s^{m}}(f^{m})_*T^-_{k-s,\varepsilon}(y)\wedge  V \Big| \leq  \\
    \int_{(\P^k)^2}& \frac{1}{d_s^{n}}-\|\psi\|_\infty \|\varphi\|_{C^2} (f^{n})^*(\omega\wedge T^+_{s,\varepsilon}(x)) \wedge \frac{1}{d_s^{m}}(f^{m})_*T^-_{k-s,\varepsilon}(y)\wedge  V, 
 \end{align*} 
 where we use that the right-hand side integral is taken outside $I(f^{n})$. Using Lemma \ref{keyCS} gives the bound:
    \begin{align*} 
 \Big| \int_{(\P^k)^2} &\chi_A(x)\xi_A(y)\psi_m(y) dd^c\varphi_n(x)\wedge  \frac{1}{d_s^{n}}(f^{n})^*T^+_{s,\varepsilon}(x) \\
  &\wedge \frac{1}{d_s^{m}}(f^{m})_*T^-_{k-s,\varepsilon}(y)\wedge  V \Big| \leq  \\
  &K \|\varphi\|_{C^2}\|\psi\|_{\infty} \sqrt{\delta}^{-n}.
 \end{align*}
 The term in $dd^c \psi_m$ is the same. That proves the point (\ref{eq1}) of the proposition. \\
 
 Assume now furtheremore that $f$ is a regular birational map. Then the proof is the same except that we have a better estimate than the one of Lemma \ref{keyCS}. Indeed, the quantity:
 $$\int_{\P^k} \frac{1}{d_{s+1}^{n}}(f^{n})^*(\omega\wedge T^+_{s,\varepsilon})\wedge U_{\frac{1}{d_s^{m}}(f^{m})_*T^-_{k-s,\varepsilon} }  $$
 converges to $\U_{T^-_{k-s}}(\frac{1}{d_{s+1}^{n}}(f^{n})^*(\omega\wedge T^+_{s}))$. The current $\frac{1}{d_{s+1}^{n}}(f^{n})^*(\omega\wedge T^+_{s})$ has mass one and support contained in $\mathrm{Supp}(T^+_s)$. Since $f$ is regular, $\U_{T^-_{k-s}}(\frac{1}{d_{s+1}^{n}}(f^{n})^*(\omega\wedge T^+_{s}))$ is uniformly bounded in $n$. The inequality (\ref{eq2}) follows. \hfill $\Box$ \hfill \\
 
  We now prove the exponential decay of correlations. \\
 
\noindent \emph{Proof of Theorem \ref{decay}.} Using  the interpolation's argument of Dinh and Sibony, it is sufficient to consider the case where $\alpha=2$ (\cite{DS4}[Corollary 6.2]).   Let $n_0:= [N/k ]$ be the integer part of $N/k$. We write $n=(k-s)n_0$ and $m=sn_0+r$ so that $0\leq r <k$ and $n+m =N$. In what follows, the constant $C_i$ changes from line to line still remaining independent of $\varphi$,  $\psi$ and $N$. Applying the previous proposition to such $n$ and $m$ implies:
\begin{align*}
\Big| \mu(\varphi\circ f^N . \psi)- \mu(\varphi)\mu(\psi) \Big| & = \Big| \mu(\varphi\circ f^n . \psi\circ f^{-m})- \mu(\varphi)\mu(\psi) \Big| \\
                                                                & \leq C\|\varphi\|_{C^2}\|\psi\|_{C^2}(\sqrt{\delta}^{-n}+\sqrt{d}^{-m}) \\
                                                                & \leq C\|\varphi\|_{C^2}\|\psi\|_{C^2}(\sqrt{\delta^{k-s}}^{-n_0}+\sqrt{d^s}^{-n_0} \sqrt{d}^{-r}).
\end{align*} 
Now, we use that $\delta^{k-s}=d^s$ and that $\sqrt{d}^{r}$ only take finitely many values. So:
\begin{align*}
\Big| \mu(\varphi\circ f^N . \psi)- \mu(\varphi)\mu(\psi) \Big| & \leq C_1\|\varphi\|_{C^2}\|\psi\|_{C^2}\sqrt{d^s}^{-n_0}\\
                                                                & \leq C_1\|\varphi\|_{C^2}\|\psi\|_{C^2}\sqrt{d}^{-s \left [\frac{N}{k}\right]}\\
                                                                &\leq  C_2\|\varphi\|_{C^2}\|\psi\|_{C^2}\sqrt{d}^{-s \left(\frac{N}{k}-1\right)}\\
                                                                &\leq  C_3\|\varphi\|_{C^2}\|\psi\|_{C^2}d^{\frac{-sN}{2k}}.
\end{align*} 
  This is exaclty what we want. The regular case is the same using the better bound of the previous proposition obtained in that case. \hfill $\Box$ \hfill \\

\noindent Gabriel Vigny, LAMFA - UMR 7352, \\ 
U. P. J. V. 33, rue Saint-Leu, 80039 Amiens, France. \\
\noindent Email: gabriel.vigny@u-picardie.fr

\end{document}